\documentclass[11pt]{article}
\usepackage{authblk}
\usepackage{amssymb}
\usepackage[all]{xy}
\usepackage{mathrsfs}
\usepackage{amsmath,amssymb,amsthm}
\usepackage{extarrows}
\usepackage[mathscr]{eucal}
\usepackage{mathrsfs}
\newtheorem{theorem}{Theorem}[section]
\newtheorem{lemma}[theorem]{Lemma}
\newtheorem{proposition}[theorem]{Proposition}
\newtheorem{definition}[theorem]{Definition}

\DeclareMathOperator{\Ima}{Im}
\DeclareMathOperator{\LIMA}{LIM}
\def\bB{\mathbf{B}}
\def\bC{\mathbf{C}}
\def\H{\mathcal{H}}
\def\C{\mathcal{C}}
\def\B{\mathcal{B}}
\def\Z{\mathcal{Z}}
\def\to{\rightarrow}
\def\d{\delta}
\def\s{\scriptstyle}
\date{}
\begin{document}
\baselineskip16pt
\title{Continuous Bounded Cohomology of Topological Semigroups}
\author{Maysam Maysami Sadr\thanks{sadr@iasbs.ac.ir}}
\affil{Department of Mathematics, Institute for Advanced Studies in Basic Sciences, Zanjan, Iran}
\maketitle
\begin{abstract}
In this short note, we give some new results on continuous bounded cohomology groups of topological semigroups with values in complex field. We show that the second continuous bounded cohomology group of a compact metrizable semigroup, is a Banach space. Also, we study cohomology groups of amenable topological semigroups, and we show that cohomology groups of rank greater than one of a compact left or right amenable semigroup, are trivial. Also, we give some examples and applications about topological lattices.

\textbf{MSC 2010.}
20J06, 22A99, 20J05, 22A26.

\textbf{Keywords.}
Topological semigroup, bounded cohomology, Banach homology.
\end{abstract}
\section{Introduction}
Homology theory is one of the most powerful tools for study of various mathematical objects.
There are many kind of (co)homology theories,
for instance: de Rham (co)homology for smooth manifolds, sheaf cohomology for algebraic varieties,
singular (co)homology of topological spaces with values in an arbitrary ring,
\u{C}ech cohomology for topological spaces, bordism homology, and Hochschild (co)homology for rings and topological algebras with values
in bimodules.

The bounded cohomology theory was first defined for discrete groups by B. E. Johnson \cite{Johns} and F. Trauber. Then, M. Gromov \cite{Gro} extended it to topological spaces. Gromov has been proved that for every path connected manifold, the bounded cohomology group of any rank is equivalent with bounded cohomology group of fundamental group of the manifold with the same rank; for more details one can look \cite{Fri2}. The continuous cohomology theory for topological spaces and topological groups have been studied by many mathematicians in different approaches; see \cite{Bro,Mdz,Mos} and \cite{Sta}.

The bounded continuous cohomology theory for topological spaces and topological groups, generalizing both continuous cohomology and bounded cohomology theories simultaneously, has been studied by many authors such as R. Frigerio \cite{Fri} and N. Monod \cite{Mon}. Bounded cohomology of semigroups has also been considered by many mathematicians such as R. Brooks \cite{B}, R. I. Grigorchuk \cite{Gri} and N. V. Ivanov \cite{I}.

In this paper we establish a topological bounded cohomology theory for topological semigroups, using continuous bounded cocycles.

In the section 2, we define bounded continuous cohomology group of a topological semigroup. In the section 3, we show some basic properties of the bounded continuous cohomology. In the next section, we explain the bounded continuous cohomology relation with amenability. In the last section, we give some examples of it.
\section{Definition of the cohomology}

For any set $X$, $\bB(X)$ denotes the Banach space of all bounded complex ($\mathbb{C}$) valued maps on $X$ with the
uniform norm. If $X$ has a topology, then $\bC(X)\subset\bB(X)$ denotes the Banach subspace of continuous maps.
By a topological semigroup we mean a semigroup $S$ with a topology  such that the multiplication $S\times S\to S$ is
jointly continuous.

Let $S$ be a semigroup. Let $\C_b^0(S)=\mathbb{C}$, and for $n\geq1$, let $\C_b^n(S)=\bB(S^n)$.
The elements of $\C_b^n(S)$ are called {\it bounded cochains} of the semigroup $S$. Let $\d^0:\C_b^0(S)\to\C^1_b(S)$
be the zero linear map and for $n\geq1$, define the bounded linear map $\d^n:\C^n_b(S)\to\C^{n+1}_b(S)$ by
\begin{equation}\label{e1}
\begin{split}
\d^n(f)(s_1,\ldots,s_{n+1})&=f(s_2,\ldots,s_{n+1})\\
&+\sum_{i=1}^n(-1)^if(s_1,\ldots,s_is_{i+1},\ldots,s_{n+1})\\
&+(-1)^{n+1}f(s_1,\ldots,s_n),
\end{split}
\end{equation}
for $f\in\C^n_b(S)$ and $s_1,\ldots,s_{n+1}\in S$. The linear map $\d^n$ is called {\it coboundary}.
It is easily checked that $\d^{n+1}\d^n=0$ and thus, we have the following cochain complex of Banach spaces and
bounded linear maps:
\begin{equation}\label{e2}
\xymatrix{0\ar[r]& \C_{b}^0(S)\ar[r]^-{{\s{\d^0}}}&\C_{b}^1(S)\ar[r]^-{{\s{\d^1}}}&
	\cdots\ar[r]^{{\s{\d^{n-1}}}}&\C_{b}^n(S)\ar[r]^-{{\s{\d^n}}}& \cdots}
	\end{equation}
Then, the cohomology groups of the complex (\ref{e2}) are called {\it bounded cohomology groups} of $S$
and denoted by $\H^n_b(S)=\frac{\ker\d^n}{\Ima\d^{n-1}}$. Always, the quotient vector space $\H^n_{cb}(S)$ is considered
as a semi normed space with quotient semi norm.
As any cohomology theory, we let $\B^n_b(S)=\Ima\d^{n-1}$ and $\Z^n_b(S)=\ker\d^n$.
The elements of $\B^n_b$ and $\Z^n_b$ are called bounded $n$-coboundaries and bounded $n$-cocycles, respectively.
For more details on bounded cohomology of semigroups, see \cite{Gri}.

Now, suppose that $S$ is a topological semigroup. Let $\C_{cb}^0(S)=\mathbb{C}$, and for $n\geq1$, let $\C_{cb}^n(S)=\bC(S^n)$. The elements of $\C_{cb}^n(S)$ are called {\it continuous bounded cochains} of
the semigroup $S$. Then we have the following Banach subcomplex of (\ref{e2}):
\begin{equation}\label{c3}
\xymatrix{0\ar[r]& \C_{cb}^0(S)\ar[r]^-{{\s{\d^0}}}&\C_{cb}^1(S)\ar[r]^-{{\s{\d^1}}}&
	\cdots\ar[r]^{{\s{\d^{n-1}}}}&\C_{cb}^n(S)\ar[r]^-{{\s{\d^n}}}& \cdots}
\end{equation}
\begin{definition}
The cohomology groups of the complex (\ref{c3}) are called  continuous bounded cohomology groups of $S$
and denoted by $\H^n_{cb}(S)$.
\end{definition}
Analogously, we have the space of continuous bounded $n$-coboundaries $\B^n_{cb}(S)$, and the space of
continuous bounded $n$-cocycles $\Z^n_{cb}$, and $\H^n_{cb}(S)$ is considered by the quotient semi norm.

{\bf Remark.}
\begin{enumerate}
\item [(I)] Let $S$ be a discrete semigroup. Consider the convolution Banach algebra $\ell^1(S)$. Then
the space $\mathbb{C}$ is a Banach $\ell^1(S)$-bimodule by the symmetric action $f\cdot\lambda=\lambda\cdot
f=\lambda\sum_{s\in S}f(s)$ for $f\in\ell^1(S)$ and $\lambda\in\mathbb{C}$. It is well known and easily
checked that the bounded Hochschild cohomology groups of $\ell^1(S)$ with values in the bimodule $\mathbb{C}$
and the bounded cohomology groups of $S$ are isometric isomorph. Thus, the bounded cohomology is a special
case of Hochschild cohomology, see\cite{M}.
\item [(II)] Let $S$ be a compact Hausdorff semigroup. If we dualize cochain complex (\ref{c3}), then (by
the natural isomorphism between $\bC(X)^*$ and the Banach space of complex Borel regular measures $\mathbf{M}(X)$
for any compact Hausdorff space $X$) we have the chain complex
$$ \xymatrix{0& \mathbb{C}\ar[l] &\mathbf{M}(S)\ar[l]_{(\s{\d^0})^*}&
	\mathbf{M}(S^2)\ar[l]_{(\s{\d^1})^*}&\cdots\ar[l]_-{(\s{\d^{2}})^*}} $$

One can consider the homology of this complex as a measure homology theory (cf. \cite{L1,L2}) that is a topological version of $\ell^1$-homology of discrete semigroups \cite{Gri}.
\end{enumerate}
\section{Some basic properties}
\begin{theorem}
Let $S,T$ be topological semigroups and $\phi:S\to T$ be a continuous homomorphism. Then for every $n\geq0$,
there is a canonical continuous linear map
$$\H^n_{cb}(\phi):\H^n_{cb}(T)\to\H^n_{cb}(S).$$
\end{theorem}
\begin{proof}
For every $n\geq1$ let $\hat{\phi}_n:\C^n_{cb}(T)\to\C^n_{cb}(S)$ be defined by
$$\hat{\phi}_n(f)(s_1,\ldots,s_n)=f(\phi(s_1),\ldots,\phi(s_n))\quad(f\in\C^n_{cb}(T)).$$
Then $(\hat{\phi}_n)_n$ is a cochain map between continuous bounded cohomology complexes of $T$ and $S$, i.e.
the following diagram is commutative:
\begin{equation}\label{e3}
\xymatrix {0\ar[r]&\mathbb{C}\ar[r]^-{\s{\d^0}}\ar[d]^{\mathrm{id}}&\C_{cb}^1(T)\ar[r]^{\s{\d^1}}\ar[d]^{\hat{\phi}_1}&
	\cdots\ar[r]^{\s{\d^{n-1}}}&\C_{cb}^n(T)\ar[r]^{\s{\d^{n}}}\ar[d]^{\hat{\phi}_n} &\cdots\\
	0\ar[r]&\mathbb{C}\ar[r]_-{\s{\d^0}}&\C_{cb}^1(S)\ar[r]_{\s{\d^1}}& \cdots \ar[r]_{\s{\d^{n-1}}} &\C_{cb}^n(S)\ar[r]_{\s{\d^{n}}}&\cdots}
\end{equation}
Thus, the standard arguments of Banach homology (\cite{H}, \cite{R}) shows that we have a continuous linear
map $\H^n_{cb}(\phi)$, defined by
$$\H^n_{cb}(\phi)(f+\B^n_{cb}(T))=\hat{\phi}_n(f)+\B^n_{cb}(S),$$
for $f\in\Z^n_{cb}(T)$.
\end{proof}
Let $\mathcal{T}\mathcal{S}\mathcal{G}$ be the category of topological semigroups and continuous homomorphisms.
Then, the above theorem shows that $\H^n_{cb}$ is a contravariant functor from $\mathcal{T}\mathcal{S}\mathcal{G}$
to the category of seminormed spaces and continuous linear maps. Since the category
$\mathcal{T}\mathcal{S}\mathcal{G}$ has no additive properties, the computation of continuous bounded cohomology
groups often are very hard. In another paper, we will consider various extensions of $\H^n_{cb}$ to some categories of
representations of topological semigroups on topological vector spaces.

For any topological semigroup $S$ it is trivial that $\H^0_{cb}(S)=\mathbb{C}$. First order cohomology groups are zero:
\begin{theorem}\label{t2}
For any topological semigroup $S$, $\H^1_{cb}(S)$ is zero.
\end{theorem}
\begin{proof}
Let $f\in\Z^1_{cb}(S)$ be a $1$-cocycle. Then for every $s,t\in S$, we have $\d^1(f)(s,t)=f(t)-f(st)+f(s)=0$ and thus,
$$f(st)=f(s)+f(t).$$
In particular, for every $s\in S$ and $n\in\mathbb{N}$, we have $f(s^n)=nf(s)$. This implies that $f(s)=0$,
since $f$ is a bounded map. Therefore $\Z^1_{cb}(S)=0$ and $\H^1_{cb}(S)$ is zero.
\end{proof}
We recall a kind of limiting process: Let $E$ be the Banach space of all
bounded sequences of complex numbers with uniform norm and let $F\subset E$
be the subspace of all convergent sequences. Then, the functional $lim:F\to\mathbb{C}$ defined by
$lim(a_n)_{n\in\mathbb{N}}=\lim_{n\to\infty}a_n$ is a bounded functional and thus, by the Hann-Banach theorem
there is a bounded functional $\LIMA:E\to\mathbb{C}$ that extends $lim$ and $\|\LIMA\|=1$ (such functionals are
called {\it Banach limits}).

\begin{theorem}\label{t1}
Let $S$ be a compact semigroup with a metric $d$ that induces the topology of $S$ and has the following property:

\begin{itemize}
\item  For every $\beta>0$,
$s,t\in S$, and $i\in\mathbb{N}$ if $d(s,t)<\beta$ then $d(s^i,t^i)<\beta$.
\end{itemize}
Then $\H^2_{cb}(S)$ is a Banach space.
\end{theorem}
\begin{proof}
It is enough that we prove $\d^1(\C^1_{cb}(S))$ is closed in $\C^2_{cb}(S)$, and thus,
it is sufficient to construct a bounded linear map $\gamma:\C^2_{cb}(S)\to\C^1_{cb}(S)$ such that
$\gamma\d=id_{\C^1_{cb}(S)}$.

Let $f\in\C^2_{cb}(S)$ be a 2-cochain. For every $s\in S$, consider the bounded sequence
$a_n^{f,s}=n^{-1}\sum_{i=1}^{n-1}f(s^i,s)$ ($n\geq2$) of complex numbers and define $\gamma(f)(s)=\LIMA(a_n^{f,s})$.
Let $\alpha>0$ be arbitrary. Since $S^2$ is a compact metric space and $f$ is continuous, there is
$\beta>0$ such that if $d(t_1,t_2)<\beta $ and $d(t_1',t_2')<\beta$ then $|f(t_1,t_1')-f(t_2,t_2')|<\alpha$.
This property together with $(*)$ implies that for every $s,t\in S$ and $n\in\mathbb{N}$ if $d(s,t)<\beta$ then
$|a_n^{f,s}-a_n^{f,t}|<\alpha$ and thus, $|\gamma(f)(s)-\gamma(f)(t)|<\alpha$. Therefore we have proved $\gamma(f)$
is continuous and $\gamma(f)\in\C^1_{cb}(S)$. Also, it is easily checked that $\gamma$ is a bounded linear operator.

Now, suppose that $g$ is in $\C^1_{cb}(S)$. For every $s\in S$ and $i\geq1$ we have
$$\d^1(g)(s^i,s)=g(s^i)-g(s^{i+1})+g(s),$$
thus, for every $n\geq2$, $a_n^{\d^1(g),s}=g(s)-n^{-1}g(s^n)$. Therefore we have
$$\gamma(\d^1(g))(s)=\LIMA(a_n^{\d^1(g),s})=\lim_{n\to\infty}g(s)-n^{-1}g(s^n)=g(s).$$
Thus, we have proved $\gamma\d(g)=g$.
\end{proof}
It is easily checked that the arguments of the proof of Theorem \ref{t1}, satisfy when $S$ is a discrete semigroup:
\begin{theorem}
Let $S$ be a discrete semigroup. Then $\H^2_{b}(S)$ is a Banach space.
\end{theorem}

\begin{proposition}\label{p1}
Let $S$ and $T$ be topological semigroups and $p:S\times T\to S$ be the natural projection. Suppose that $T$
has a unite element $e$. Then the linear map $\H^n_{cb}(p):\H^n_{cb}(S)\to\H^n_{cb}(S\times T)$
is injective for all $n\geq1$.
\end{proposition}
\begin{proof}
For $n=1$, the result follows from Theorem \ref{t2}. Let $n\geq2$ be fixed, and let $p^{(n)}:(S\times T)^n\to S^n$
be defined by
$$p^{(n)}((s_1,t_1),\ldots,(s_n,t_n))=(s_1,\ldots,s_n)$$
for $s_1,\ldots,s_n\in S, t_1,\ldots,t_n\in T$.
By definition of $\H^n_{cb}(p)$, we must prove that if $f\in\Z^n_{cb}(S)$ and $f\circ p^{(n)}\in\B^n_{cb}(S\times T)$,
then $f$ is in $\B^n_{cb}(S)$. Thus, consider such a $n$-cocycle $f$. There is $g\in\C^{n-1}_{cb}(S\times T)$
such that $\d^{n-1}(g)=f\circ p^{(n)}$. Define $\hat{g}\in\C^{n-1}_{cb}(S)$ by
$$\hat{g}(s_1,\ldots,s_{n-1})=g((s_1,e),\ldots,(s_{n-1},e))\quad\quad\quad(s_1,\ldots,s_{n-1}\in S).$$
Then, for every $s_1,\ldots,s_n\in S$, we have
\begin{equation*}
\begin{split}
\d^{n-1}(\hat{g})(s_1,\ldots,s_{n})&=\hat{g}(s_2,\ldots,s_{n})
+\sum_{i=1}^{n-1}(-1)^i\hat{g}(s_1,\ldots,s_is_{i+1},\ldots,s_{n})\\
&+(-1)^{n}\hat{g}(s_1,\ldots,s_{n-1})
=g((s_2,e),\ldots,(s_{n},e))\\
&+\sum_{i=1}^{n-1}(-1)^ig((s_1,e),\ldots,(s_is_{i+1},e),\ldots,(s_{n},e))\\
&+(-1)^{n}g((s_1,e),\ldots,(s_{n-1},e))\\
&=\d^{n-1}(g)((s_1,e),\cdots,(s_{n},e)).
\end{split}
\end{equation*}
On the other hand, $\d^{n-1}(g)((s_1,e)\ldots,(s_{n},e))=f(s_1,\ldots,s_{n})$.
Thus, we have $\d^{n-1}(\hat{g})=f$ and $f\in\B^n_{cb}(S)$.
\end{proof}
\section{Relation with amenability}
Let $S$ be a topological semigroup. A function $f\in\bC(S)$ is called {\it right uniformly continuous}, if the map
$\Phi_f:S\to\bC(S)$ defined by $\Phi_f(s)=f\cdot s$ is continuous with uniform norm of $\bC(S)$, where
$f\cdot s(x)=f(sx)$ ($x\in S$). {\it Left} uniformly continuous functions are similarly defined.
The space of all right (left)
uniformly continuous functions is denoted by $\mathbf{R}\mathbf{U}\bC(S)$ ($\mathbf{L}\mathbf{U}\bC(S)$). Note that
if $f\in\mathbf{R}\mathbf{U}\bC(S)$ and $s\in S$, then $f\cdot s\in\mathbf{R}\mathbf{U}\bC(S)$. Also, it is easily
checked that $\mathbf{R}\mathbf{U}\bC(S)=\mathbf{L}\mathbf{U}\bC(S)=\bC(S)$ when $S$ is compact, and it
is clear that $\mathbf{R}\mathbf{U}\bC(S)=\mathbf{L}\mathbf{U}\bC(S)=\bC(S)=\bB(S)$ when $S$ is discrete.

A topological semigroup $S$ is called {it left amenable} if there is a {\it left invariant mean} on $\mathbf{R}\mathbf{U}\bC(S)$, i.e. a bounded linear functional $m$ on $\mathbf{R}\mathbf{U}\bC(S)$ such that
$\langle m,1_S\rangle=\|m\|=1$ (where $1_S$ is the constant map on $S$ with value $1$) and for every $s\in S$ and
$f\in\mathbf{R}\mathbf{U}\bC(S)$, $\langle m,f\cdot s\rangle=\langle m,f\rangle$. {\it Right} invariant means and
{\it right} amenable semigroups are similarly defined. A topological semigroup is called amenable if it is both left
and right amenable.

It is well known and easily checked that for topological semigroups $S$ and $T$, if there is a continuous
homomorphisms form $S$ onto $T$, and $S$ is left (right) amenable, then $T$ is also left (right) amenable.
In particular, if $S$ is left (right) amenable semigroup with topology $\tau$, and $\tau'$ is another
semigroup topology on $S$ such that $\tau'\subset \tau$, then $(S,\tau')$ is left (right) amenable.
Thus, any commutative topological semigroup is amenable since any commutative discrete semigroup is amenable
(\cite{P}). It is well known that any compact group is amenable (\cite{P}), but there are compact semigroups
that are not left amenable nor right amenable:

{\bf Example.} Let $X$ and $Y$ be two disjoint compact spaces with distinguished elements $x_0\in X$
and $y_0\in Y$. Define a semigroup multiplication on disjoint union space $T=X\cup Y$ by
$$xx'=x_0,\quad yy'=y_0,\quad xy=x_0,\quad yx=y_0,$$
for every $x,x'\in X$ and $y,y'\in Y$. Then $T$ becomes a compact semigroup. We show that $T$ is not left amenable.
Suppose $m$ is a bounded linear functional on $\bC(T)$ such that $\langle m,1_T\rangle=\|m\|=1$. For every
$x\in X\subset T$ and $f\in\bC(T)$, we have $\langle m,f\cdot x\rangle=\langle m,f(x_0)1_T\rangle=f(x_0)$, and
similarly $\langle m,f\cdot y\rangle=f(y_0)$ for every $y\in Y$. Thus, $m$ is not a left invariant mean,
since there is a continuous map $f$ on $S$ such that $f(x_0)\neq f(y_0)$. Thus, we have proved that $T$ is not
left amenable. Let $T^{op}$ be the opposite semigroup of $T$. Then $T^{op}$ is not right amenable. Now
the compact semigroup $S=T\times T^{op}$ is not left nor right amenable, since the canonical projection maps from $S$
to $T$ and $T^{op}$ are continuous surjective homomorphisms.

We need the following simple topological lemma.
\begin{lemma}\label{l1}
Let $X$ be a topological space and $Y$ be a compact space. Let $f:X\times Y\to\mathbb{C}$ be a continuous map.
Then $F:X\to\bC(Y)$, defined by $F(x)(y)=f(x,y)$ is  continuous with norm topology of $\bC(Y)$.
\end{lemma}
\begin{proof}
Let $x_0\in X$ and $\alpha>0$ be arbitrary. Since $f$ is continuous, for every $y\in Y$, there are open
sets $U_y, V_y$ in $X$ and $Y$ respectively, such that $(x_0,y)\in U_y\times V_y$ and
$|f(x_0,y)-f(x,y')|<\alpha/2$ for every $(x,y')\in U_y\times V_y$. Since $Y$ is compact, there are
$y_1,\ldots,y_n\in Y$ such that $Y=\cup_{i=1}^nV_{y_i}$. Let $W$ be the open set $\cap_{i=1}^nU_{y_i}$.
Let $x\in W$ and $y\in Y$ be arbitrary. Then for some $i$ ($i=1,\cdots,n$), $y$ belongs to $V_{y_i}$ and we have,
\begin{equation*}
\begin{split}
&|f(x,y)-f(x_0,y)|\leq\\
&|f(x,y)-f(x_0,y_i)|+|f(x_0,y_i)-f(x_0,y)|<\alpha/2+\alpha/2=\alpha.
\end{split}
\end{equation*}
Thus, we have $\|F(x)-F(x_0)\|<\alpha$ for every $x\in W$. The proof is complete.
\end{proof}
The proof of the following Theorem is an adaptation of the proof given in \cite[Theorem 2.1]{Gri}
to the topological case.
\begin{theorem}\label{t3}
Let $S$ be a compact semigroup and suppose that $S$ is left (right) amenable.
Then $\H^n_{cb}$ is zero for every $n\geq0$.
\end{theorem}
\begin{proof}
Suppose that $S$ is left amenable and let $m$ be a left invariant mean on $\bC(S)^*$. Similar \cite{Gri},
we use the notation
\begin{equation}\label{e8}
m(f)=\int_Sf(s)d(s)
\end{equation}
for $f\in\bC(S)$. Thus, we have
\begin{enumerate}
\item [(i)] $\int_S1_S(s)d(s)=1$, and
\item [(ii)] $\int_Sf(ts)d(s)=\int_Sf(s)d(s)$ for every $f\in\bC(S)$ and $t\in S$.
\end{enumerate}
The cases $n=0$ and $n=1$ were considered before, thus, suppose that $n\geq2$ and let $f\in\Z^n_{cb}(S)$. Then, for
every $s_1,\ldots,s_{n+1}\in S$, we have,
\begin{equation*}
\begin{split}
\d^n(f)(s_1,\ldots,s_{n+1})&=f(s_2,\ldots,s_{n+1})\\
&+\sum_{i=1}^n(-1)^if(s_1,\ldots,s_is_{i+1},\ldots,s_{n+1})\\
&+(-1)^{n+1}f(s_1,\ldots,s_n)=0\\
\end{split}
\end{equation*}
If we fix $s_1,\ldots,s_n\in S$ and integrate the above formula over the variable $s_{n+1}$ in the sense of
(\ref{e8}), then we have
\begin{equation}\label{e9}
\begin{split}
&\int_Sf(s_2,\ldots,s_{n+1})d(s_{n+1})\\
+&\sum_{i=1}^{n-1}(-1)^i\int_Sf(s_1,\ldots,s_is_{i+1},\ldots,s_n,s_{n+1})d(s_{n+1})\\
+&(-1)^{n}\int_Sf(s_1,\ldots,s_{n-1},s_ns_{n+1})d(s_{n+1})\\
+&(-1)^{n+1}\int_Sf(s_1,\cdots,s_n)d(s_{n+1})=0\\
\end{split}
\end{equation}
By property (i),
\begin{equation}\label{e10}
\int_Sf(s_1,\ldots,s_n)d(s_{n+1})=f(s_1,\ldots,s_n),
\end{equation}
and by property (ii),
\begin{equation}\label{e11}
\int_Sf(s_1,\ldots,s_{n-1},s_ns_{n+1})d(s_{n+1})=\int_Sf(s_1,\ldots,s_{n-1},s_{n+1})d(s_{n+1}).
\end{equation}
Let $g:S^{n-1}\to\mathbb{C}$ be defined by
$$g(s_2,\ldots,s_{n})=\int_Sf(s_2,\ldots,s_n,s_{n+1})d(s_{n+1}).$$
By Lemma \ref{l1}, the map $F:S^{n-1}\to\bC(S)$, defined by
$$F(s_2,\ldots,s_{n})(x)=f(s_2,\ldots,s_{n},x)\quad\quad(x\in S),$$
is continuous with the norm of $\bC(S)$. On the other hand, $\int:\bC(S)\to\mathbb{C}$ is also continuous
with the norm. Thus, the map $g=\int F$ is in $\C^{n-1}_{cb}(S)$. Therefore, by (\ref{e9}), (\ref{e10})
and (\ref{e11}), we have,
\begin{equation*}
\begin{split}
(-1)^{n}f(s_1,\ldots,s_{n})&=g(s_2,\ldots,s_{n})\\
&+\sum_{i=1}^{n-1}(-1)^ig(s_1,\ldots,s_is_{i+1},\ldots,s_{n})\\
&+(-1)^{n}g(s_1,\ldots,s_{n-1})
\end{split}
\end{equation*}
But the right hand side of the latter equation is $\d^{n-1}(g)$. Thus,
$$f=\d^{n-1}((-1)^ng).$$
Therefore we have proved $\B^n_{cb}(S)=\Z^n_{cb}(S)$ and $\H^n_{cb}(S)=0$. A similar proof can be given in
the case of right amenable $S$.
\end{proof}
\section{Some examples}
Gromov (\cite{Gro}) proved that for any connected manifold $X$, and any $n\geq1$, the bounded cohomology
of $X$ and the bounded cohomology of the fundamental homotopy group $\pi_1(X)$ of $X$ coincide
(for more details see \cite{Gro,B,I}, and \cite{Gri}). Thus, there are many discrete
groups that their bounded cohomology groups are non zero.

Let $G$ be a discrete group and $S$ be a topological semigroup with a unite. Suppose that for an integer $n\geq2$,
$\H^n_{cb}(G)\neq0$ (for example $G=F_2$, the free group on two generators, and $n=2$, see \cite{Gri}, \cite{F}).
Then by Proposition \ref{p1} we have $\H^n_{cb}(G\times S)\neq0$.

A semigroup $S$ is called {\it semilattice} if it is commutative and $ss=s$ for every $s\in S$.
\begin{theorem}
Let $S$ be a topological semilattice. Then $\H^2_{cb}(S)$ is zero.
\end{theorem}
\begin{proof}
Let $f\in\Z^2_{cb}(S)$ be a $2$-cocycle. We need a $g\in\C^1_{cb}(S)$ such that for every $s,t\in S$,
$$f(s,t)=g(s)+g(t)-g(s,t).$$
Since $f$ is a $2$-cocycle, for every $s_1,s_2,s_3\in S$ we have
\begin{equation}\label{e4}
\d^2(f)(s_1,s_2,s_3)=f(s_2,s_3)-f(s_1s_2,s_3)+f(s_1,s_2s_3)-f(s_1,s_2)=0.
\end{equation}
Apply (\ref{e4}) with $s_1=s,s_2=s,s_3=t$, we obtain
\begin{equation}\label{e5}
f(s,s)=f(s,st)
\end{equation}
and similarly
\begin{equation}\label{e6}
f(t,t)=f(t,st).
\end{equation}
Apply (\ref{e4}) with $s_1=s,s_2=t,s_3=st$, we obtain
\begin{equation}\label{e7}
f(s,t)=f(t,st)-f(st,st)+f(s,st)
\end{equation}
Now, by (\ref{e5}), (\ref{e6}) and (\ref{e7}), we have
$$f(s,t)=f(t,t)+f(s,s)-f(st,st).$$
Thus, if we define $g(s)=f(s,s)$ ($s\in S$), then $g\in\C^1_{cb}(S)$ and $\d^1(g)=f$.
\end{proof}
{\bf Remark.}
\begin{enumerate}
\item [(a)] The above result follows directly from Theorem \ref{t3}, when $S$ is compact.
\item [(b)] In \cite{C}, Y. Choi proved , that for any discrete semilattice $S$ and any symmetric Banach
$\ell^1(S)$ bimodule $E$, every Hochschild  cohomology group of $\ell^1(S)$ with coefficient in $E$ is trivial
(see Remark of Section 1). Thus, for every discrete semilattice $S$ and $n\geq0$, $\H^n_b(S)=0$.
\end{enumerate}

{\bf Example.} Let $(X,d)$ be a metric space and let for every $A\subset X$ and $\epsilon>0$, $N_\epsilon(A)$
be the $\epsilon$-neighborhood of $A$ in $X$. Let $S_X$ be the set of all nonempty closed bounded subsets of $X$.
Then $S_X$ is a metric space by the following metric that is called Hausdorff distance:
$$d_H(C_1,C_2)=\inf\{\epsilon>0: C_1\subset N_\epsilon(C_2)\text{ and } C_2\subset N_\epsilon(C_1)\},$$
for $(C_1,C_2\in S_X)$.
Also $S_X$ is a semilattice with the multiplicaton $C_1C_2=C_1\cup C_2$. It is easily checked that $S_X$
is a topological semilattice with the topology induced by $d_H$. Note that if $(X,d)$ is a compact metric
space then so is $(S_X,d_H)$, \cite[Lemma 5.31]{BH}.\\
Let $(X,d)$ and $(Y,\rho)$ be disjoint compact metric spaces with distinguished elements $x_0\in X$ and
$y_0\in Y$. Let $S$ and $T$ be compact semigroups defined in Example of Section 3, using $X$ and $Y$. Define
the metric $d'$ on $T=X\cup Y$ by $d'|_{X\times X}=d$, $d'|_{Y\times Y}=\rho$  and $d'(x,y)=1$ for $x\in X,y\in Y$.
Then $T$ together with $d'$, satisfy the conditions of Theorem \ref{t1}. Also, $S$ together with the maximum metric
satisfy the conditions of Theorem \ref{t1}, but it is not left nor right amenable, and thus, we can not apply Theorem
\ref{t3}, to conclude that the cohomology groups of $S$ are zero.
\bibliographystyle{amsplain}

\end{document}